\documentclass[11pt]{amsart}
\usepackage{geometry}                
\usepackage{graphicx}
\usepackage{amssymb}
\usepackage{epstopdf}
\usepackage[all]{xypic}
\usepackage{hyperref}
\usepackage[arrow,curve,matrix]{xy}

\theoremstyle{plain}

\newtheorem{thm}{Theorem}[section]

\newtheorem{lem}[thm]{Lemma}

\theoremstyle{definition}

\newtheorem{ex}[thm]{Example}
\newtheorem{rmk}[thm]{Remark}

\title{Lazarsfeld-Mukai bundles and applications. II}
\author{Marian Aprodu}
\date{\today}
\dedicatory{To Alexandru Dimca and \c Stefan Papadima}                                       
\thanks{This work was partly supported by a Humboldt Fellowship and by the CNCS/UEFISCDI grant PN-II-ID-PCE-2011-3-0288 contract no. 132/05.10.2011. I would like to thank Humboldt Universit\"at for hospitality and to G. Farkas for useful discussions on the topic.}
\address{"Simion Stoilow" Institute of Mathematics of the Romanian Academy, P.O. Box 1-764, RO-014700, Bucharest, Romania}
\email{Marian.Aprodu@imar.ro}

\begin{document}
\maketitle

\section{Introduction}
The notion of Lazarsfeld-Mukai bundle goes back to the 1980's, when two important problems in algebraic geometry were solved using vector bundle techniques \cite{Lazarsfeld86,Mukai}. They were initially defined as vector bundles with particularly nice properties on $K3$ surfaces, and their main applications to date remain within the $K3$ framework. The definition makes sense however in a much larger class of surfaces.

\begin{figure}[h!]
\label{fig:LM}
  \centering
    \includegraphics[width=0.9\textwidth]{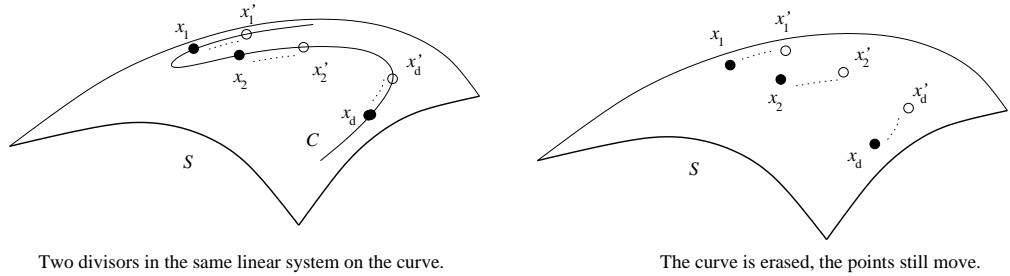}
\caption{Divisors in $|A|$ as moving subschemes of $S$.}
\end{figure}

Let $S$ be a surface, $C$ be a smooth curve on $S$, and $A$ be a $g^1_d$ on $C$. A natural question related to this setup is the following: can $A$ be lifted to the surface $S$? The chances for $A$ to be induced by a pencil on $S$ are slim, by the simple fact that we cannot exclude the possibility that $\mathrm{Pic}(S)$ be generated by $C$ itself. This case actually occurs in many situations, for instance on very general $K3$ surfaces. Instead, we can try and lift the pencil $|A|$ to a {\em different object}, and in doing so we have to decide what kind of {\em object} would that be. If we recall that an element $D$ in the linear system $|A|$ is a sum of points $x_1+\cdots+x_d$ on the curve $C$, hence points on the surface $S$, then it is plausible that any such collection of points on $S$ be cut out by a section in a rank-two vector bundle. Since the points move in a linear system on $C$, the section would also have move in a two-dimensional space of sections over $S$  in this hypothetical bundle. In other words, we want to erase the curve from the picture, keep the moving points and interpret them as elements in a linear family of zero-dimensional subschemes on $S$, see Figure~\ref{fig:LM}. If this goal is achieved, then the rank-two bundle in question, which comes equipped with a natural two-dimensional space of global sections, is called a {\em Lazarsfeld-Mukai} bundle \cite{Mukai,Lazarsfeld86,Lazarsfeld87}.
This strategy works well for regular surfaces \cite{Pareschi95}, and can be extended for higher-dimensional linear systems on $C$, \cite{Lazarsfeld86,Lazarsfeld87}. The details are explained in Section~\ref{sec:LM}.

Lazarsfeld-Mukai bundles proved to be useful in various situations. They appear in the classification of Fano threefolds \cite{Mukai}, in classical Brill-Noether problems \cite{Lazarsfeld86,Lazarsfeld87}, in higher Brill-Noether theory \cite{Farkas-Ortega}, in syzygy-related questions \cite{Voisin02,Voisin05,Aprodu-Farkas} etc; see \cite{Lazarsfeld87,Aprodu13} for some surveys of this topic.

The specific problem we consider here is the computation of the dimensions of Brill-Noether loci. For curves on regular surfaces, this computation reduces to a dimension calculation of the parameter spaces of Lazarsfeld-Mukai bundles. Beyond the Brill-Noether theoretical interest for this type of problems, the motivation comes from syzygy theory, see Section \ref{sec:syz} for a more detailed discussion. Our goal is twofold. On the one hand, we review the general theory that is generally focused on the case of $K3$ surfaces. This recollection of facts might be of future use. On the other hand, we slightly improve some results obtained so far in the non-$K3$ case.

The present work is a continuation of \cite{Aprodu13} and its outline is the following. In Section~\ref{sec:LM} we recall the definition of Lazarsfeld-Mukai bundles on regular surfaces and we prove a general dimension statement, Theorem~\ref{thm:rho}: the Brill-Noether loci have the expected dimension if some suitable vanishing conditions are satisfied. In Section~\ref{sec:LM rank 2} we estimate the dimensions of Brill-Noether loci for curves on rational surfaces with an anticanonical bound, Theorem \ref{thm:main}. It is an extension of the main result of \cite{LelliChiesa}. In Section~\ref{sec:syz} we discuss some applications to syzygies, based on the main result of \cite{Aprodu05}, Theorem~\ref{thm:Green}. An alternate proof of  \cite[Theorem 8.1]{Aprodu02} is given in Example~\ref{ex:Hirzebruch}.

\section{Lazarsfeld-Mukai bundles on surfaces with $q=0$}
\label{sec:LM}

We follow closely the approach from \cite{Pareschi95,Lazarsfeld87}. Let $S$ be a surface with $h^1(\mathcal O_S)=0$, $C$ be a smooth connected curve of genus $g$ in $S$ and denote $L=\mathcal O_S(C)$. The hypothesis $h^1(\mathcal O_S)=0$ is needed for technical reasons and ensures that $T_C|L|=H^0(N_{C|S})$.
Let $A$ be a base-point-free complete $g^r_d$ on $C$ and denote by $M_A$ the kernel of
\[
\mathrm{ev}_A:H^0(A)\otimes \mathcal O_C\to A.
\] 

The evaluation map lifts to a surjective sheaf morphism $H^0(A)\otimes \mathcal O_S\to A$ on $S$
whose kernel $F_{C,A}$ is a vector bundle of rank $(r+1)$. Its dual $E_{C,A}=F_{C,A}^*$ is called a {\em Lazarsfeld-Mukai bundle}. Dualizing the defining sequence of $E_{C,A}$
\begin{equation}
\label{eqn:F}
0\to F_{C,A}\to H^0(A)\otimes \mathcal O_S\to A\to 0
\end{equation}
we obtain the defining sequence of $E_{C,A}$
\begin{equation}
\label{eqn:E}
0\to H^0(A)^*\otimes \mathcal O_S\to E_{C,A}\to N_{C|S}\otimes A^*\to 0.
\end{equation}

The following properties of $E_{C,A}$ and $F_{C,A}$ are obtained by direct computation, using the hypotheses $h^1(\mathcal{O}_S)=0$ and $h^0(C,A)=r+1$ \cite{Lazarsfeld87}, \cite{Pareschi95}, \cite{LelliChiesa}:
\begin{enumerate}
\item
$\mathrm{det}(E_{C,A})=L$,

\item
$c_2(E_{C,A})=d$,

\item
$h^0(S,F_{C,A})=h^1(S,F_{C,A})=0$, 

\item
$\chi(S,F_{C,A})=h^2(S,F_{C,A})=(r+1)\chi(\mathcal O_S)+g-d-1$,

\item
$h^0(S,E_{C,A})=r+1+h^0(C,N_{C|S}\otimes A^*)$,

\item
$E_{C,A}$ is generated off the base locus of $|N_{C|S}\otimes A^*|$ inside $C$.
\end{enumerate}

Restricting the sequence (\ref{eqn:F}) to the curve $C$, we obtain a short exact sequence:
\begin{equation}
\label{eqn:F|C}
0\to N_{C|S}^*\otimes A\to F_{C,A}|_C\to M_A\to 0
\end{equation}
which implies, twisting by $K_C\otimes A^*$ and using the adjunction formula,
\begin{equation}
\label{eqn:F|C al doilea}
0\to \omega_S|_C\to F_{C,A}\otimes K_C\otimes A^*\to M_A\otimes K_C\otimes A^*\to 0.
\end{equation}

Note that $H^0(M_A\otimes K_C\otimes A^*)=\ker(\mu_{0,A})$, where $\mu_{0,A}:H^0(A)\otimes H^0(K_C\otimes A^*)\to H^0(K_C)$ is the Petri map.

Recall \cite{Arbarello-Cornalba_Petri} that, for any $r$ and $d$, the Brill--Noether loci form a family $\mathcal W^r_d(|L|)$ over the open subset of $|L|$ corresponding to smooth curves.
The next result uses the hypothesis $h^1(\mathcal{O}_S)=0$ and follows from the discussion in \cite[p. 197]{Pareschi95}  (see also \cite[Lemma 2.3]{Aprodu-Farkas}):

\begin{lem}
\label{lem:Pareschi}
If $(C,A)\in \mathcal W$ is a general pair in an irreducible component  of $\mathcal W^r_d(|L|)$ dominating over $|L|$, then the coboundary map $H^0(C,M_A\otimes K_C\otimes A^*)\to H^1(C,\omega_S|_C)$ vanishes.
\end{lem}

Lemma \ref{lem:Pareschi} exhibits an exact sequence
\begin{equation}
\label{eqn:generic}
0\to H^0(C,\omega_S|_C)\to H^0(C,F_{C,A}\otimes K_C\otimes A^*)\to \ker(\mu_{0,A})\to 0
\end{equation}
for a general choice of a pair $(C,A)\in \mathcal W$. In particular, $\mathcal W$ is a smooth of expected dimension at the point $(C,A)$ if the following equality holds: $h^0(C,\omega_S|_C)=h^0(C,F_{C,A}\otimes K_C\otimes A^*)$. The sequence (\ref{eqn:generic}) is useful to estimate the dimension of Brill-Noether loci and, in some situations, smoothness follows from appropriate vanishing conditions:

\begin{thm}
\label{thm:rho}
Notation as above. Assume that $h^2(S,L)=0$. Let $(C,A)$ be a general pair in a dominating component $\mathcal W$ such that $h^2(S,F_{C,A}\otimes E_{C,A})=h^2(S,\mathcal O_S)$ and $h^2(S,E_{C,A})=0$. Then $\mathcal W$ is of dimension $\le \rho(g,r,d)+\mathrm{dim}|L|+(r+1)h^1(S,E_{C,A})$ at the point $(C,A)$. In particular, if $h^1(S,E_{C,A})=0$, then $\mathcal W$ is smooth of expected dimension $\rho(g,r,d)+\mathrm{dim}|L|$ at the point~$(C,A)$. 
\end{thm}

\proof
From the long exact sequence associated to the sequence
\[
0\to \mathcal O_S\to L\to N_{C|S}\to 0,
\]
applying the vanishing hypothesis $h^2(S,L)=0$ we obtain an a surjection $H^1(C,N_{C|S})\to H^2(S,\mathcal O_S)$ and hence, since $K_C\cong N_{C|S}\otimes \omega_S|_C$, we have $h^0(C,\omega_S|_C)\ge h^2(S,\mathcal O_S)$. From the sequence (\ref{eqn:generic}) it follows that 
$\mathrm{dim}(\ker(\mu_{0,A}))\le h^0(C,F_{C,A}\otimes K_C\otimes A^*)-h^2(\mathcal O_S)$.

Twisting the sequence (\ref{eqn:E}) by $F_{C,A}\otimes \omega_S$, taking global sections in the sequence
\begin{equation}
\label{eqn:ker mu}
0\to H^0(A)^*\otimes F_{C,A}\otimes \omega_S\to F_{C,A}\otimes E_{C,A}\otimes \omega_S\to F_{C,A}\otimes K_C\otimes A^*\to 0
\end{equation}
applying Serre duality and using the hypothesis: $h^0(F_{C,A}\otimes \omega_S)=0$ and $h^0(F_{C,A}\otimes E_{C,A}\otimes \omega_S)=h^2(\mathcal O_S)$ we obtain the inequality $h^0(C,F_{C,A}\otimes K_C\otimes A^*)\le h^2(\mathcal O_S)+(r+1)h^1(S,E_{C,A})$. 

We obtain $\mathrm{dim}\ker(\mu_{0,A})\le (r+1)h^1(S,E_{C,A})$ and hence the conclusion follows.
\endproof

\begin{rmk}
\label{rmk:ad}
The assumption $h^2(S,F_{C,A}\otimes E_{C,A})=h^2(\mathcal O_S)$ is natural. For stable bundles, it is a sufficient condition for the smoothness of the moduli space at the point defined by $E_{C,A}$. Absent this condition, we obtain the weaker estimate
\[
\mathrm{dim}_{(C,A)}\mathcal W\le \rho(g,r,d)+\mathrm{dim}|L|+(r+1)h^1(S,E_{C,A})+(h^2(S,F_{C,A}\otimes E_{C,A})-h^2(\mathcal O_S)).
\]
Note that $\mathcal O_S$ is a direct summand of $F_{C,A}\otimes E_{C,A}$, its complement is $\mathrm{ad}(E_{C,A})$, the bundle of trace-free endomorphisms, and hence 
$$h^2(S,F_{C,A}\otimes E_{C,A})-h^2(\mathcal O_S)=h^2(S,\mathrm{ad}(E_{C,A})).$$
\end{rmk}

\begin{rmk}
\label{rmk:pg=0}
If in addition $p_g(S)=h^2(\mathcal O_S)=0$ then the vanishing of $h^2(S,E_{C,A})$ follows from $h^2(S,F_{C,A}\otimes E_{C,A})=0$ and Serre duality in the sequence (\ref{eqn:ker mu}). The condition $h^2(S,L)=h^0(S,L^*\otimes \omega_S)=0$ is also automatic, as $L$ is effective and hence $h^0(S,L^*\otimes \omega_S)\le h^0(S,\omega_S)$. Furthermore, the sequence (\ref{eqn:F}) twisted by $\omega_S$ shows that $h^1(S,E_{C,A})=h^0(C,\omega_S\otimes A)$ in this case.
\end{rmk}

The hypotheses of Theorem \ref{thm:rho} are realised in a number of situations, which we enumerate below.

\begin{enumerate}
\item
{\em $K3$ surfaces}, \cite{Lazarsfeld86,Lazarsfeld87}, see also \cite{Pareschi95,Aprodu-Farkas}. In this case, $\omega_S\cong\mathcal O_S$ and $N_{C|S}=K_C$. The hypothesis $h^2(S,F_{C,A}\otimes E_{C,A})=h^2(\mathcal O_S)=1$ is equivalent to simplicity of $E_{C,A}$. The vanishing of $h^1(S,E_{C,A})$ and of $h^2(S,E_{C,A})$ follow from Serre duality and the vanishing of $h^1(S,F_{C,A})$ and of $h^0(S,F_{C,A})$.

\item
{\em Enriques surfaces}, compare to \cite{Rasmussen-Zhou}. Assume $S$ is an Enriques surface and consider $X\to S$ the $K3$ universal cover of $S$. Suppose that for general pair $(C,A)$ in a dominating component $\mathcal W$,  the associated Lazarsfeld-Mukai bundle $E_{C,A}$ is stable with respect to a given polarization $H$. Since the property of being a Lazarsfeld-Mukai bundle is an open condition, the main result of \cite{Kim98} (see Theorem on page 88) shows that for a general $(C,A)$, the bundle $E_{C,A}$ is not isomorphic to $E_{C,A}\otimes\omega_S$ and hence $h^2(S,F_{C,A}\otimes E_{C,A})=h^2(\mathcal O_S)=0$. As noted in Remark \ref{rmk:pg=0}, the condition $h^2(S,E_{C,A})=0$ follows.

\item
{\em Rational surface with an anti-canonical pencil}, \cite{LelliChiesa}. Suppose that $E_{C,A}$ is stable with respect to a given polarization $H$ and $A\not\cong\omega_S^*|_C$. Since $\omega_S^*$ is effective, it follows that there are no non-zero morphisms from $E_{C,A}\otimes \omega_S^*$ to $E_{C,A}$, hence $h^2(S,F_{C,A}\otimes E_{C,A})=h^2(\mathcal O_S)=0$. As pointed out in \cite{LelliChiesa}, the existence of an anti-canonical pencil implies $h^1(S,E_{C,A})=0$. The vanishing of $h^2(S,E_{C,A})=h^0(S,F_{C,A}\otimes \omega_S)$ follows from the vanishing of $h^0(S,F_{C,A})$.
\end{enumerate}

\section{Lazarsfeld-Mukai bundles of rank two on rational surfaces}
\label{sec:LM rank 2}

It was pointed out in Remark~\ref{rmk:pg=0} that if $q(S)=p_g(S)=0$, some of the hypotheses of Theorem~\ref{thm:rho} are superfluous. Under this assumption, since there are fewer conditions, Theorem~\ref{thm:rho} can be substantially refined. We shall consider the case of rank-two Lazarsfeld-Mukai bundles, i.e. associated to base-point-free complete $g^1_d$'s, on  rational surfaces. The main result is the following:

\begin{thm}
\label{thm:main}
Let $S$ be a rational surface, $L\ge 0$ be an effective line bundle on $S$ and denote by $k\ge 3$ the maximal gonality of smooth curves on $S$ and by $g$ their genus. Let $H$ be an ample line bundle on $S$ and $C$ be a smooth curve on $S$ of gonality $k$. Suppose that $\omega_S\cdot H\le 0$, $-\omega_S\cdot C\ge k$ and that $C$ is of Clifford dimension one. Then, for any integer $d$ such that $k\le d\le g-2k+2$, and any component $\mathcal W$ of $\mathcal W^1_d(|L|)$ that dominates the linear system $|L|$, we have
\begin{equation}
\label{eqn:lgc W}
\mathrm{dim}(\mathcal W)\le \mathrm{dim}|L|+(d-k).
\end{equation}
\end{thm}

Note that of $C$ is an ample divisor, then the polarisation $H$ can be chosen to be $\mathcal{O}_S(C)$. If  the anticanonical bundle is effective,  then the polarisation $H$ can be arbitrarily chosen.

Lelli-Chiesa proved this result for rational surfaces $S$ with $h^0(S,\omega_S^*)\ge 2$. Note that this condition implies automatically $-\omega_S\cdot C\ge k$, as an anticanonical pencil will restrict to a pencil on $C$. For sake of completeness, we present here a full proof covering also the case of surfaces with an anticanonical pencil \cite{LelliChiesa}, using the strategy from \cite{Aprodu-Farkas}; the really new case compared to \cite{LelliChiesa} is Subcase I.c.

\proof
We proceed by induction on $d\ge k$. There are several possible cases, according to the behaviour of the $g^1_d$'s and of their associated Lazarsfeld-Mukai bundles.

\medskip

{\em Case I.} 
Assume that, for  $(C,A)\in \mathcal W$ general, $A$ is base-point-free and complete. 

\medskip

{\em Subcase I.a.} 
Assume that, for $(C,A)\in \mathcal W$ general, $h^1(S,E_{C,A})=0$ and $E_{C,A}$ is $H$--stable. Since $\omega_S^*\cdot H\ge 0$, the stability implies that there is no non-zero morphism form $E_{C,A}\otimes \omega_S^*$ to $E_{C,A}$, i.e. $h^2(S,F_{C,A}\otimes E_{C,A})=h^0(S,F_{C,A}\otimes E_{C,A}\otimes \omega_S)=0$. We apply Remark \ref{rmk:pg=0} and Theorem \ref{thm:rho}.

\medskip

{\em Subcase I.b.} 
Assume that, for $(C,A)\in \mathcal W$ general, $h^1(S,E_{C,A})=0$ and $E_{C,A}$ is not $H$--stable. If $h^2(S,F_{C,A}\otimes E_{C,A})=0$, we apply again Theorem \ref{thm:rho}, hence we may assume than there is a non-zero morphism from $E_{C,A}\otimes \omega_S$ to $E_{C,A}$. The bundle $E_{C,A}$ has a maximal destabilising subsheaf $M$, which induces an extension
\begin{equation}
\label{eqn:max dest}
0\to M\to E_{C,A}\to N\otimes\mathcal I_\xi\to 0,
\end{equation}
with $M\cdot H\ge N\cdot H$, where $\xi$ is a zero-dimensional locally complete intersection subscheme of length $\ell=\ell(\xi)=M\cdot N-d$. Note that 
if $E_{C,A}\not\cong M\oplus N$ then we have (compare to \cite{Aprodu-Farkas} Lemma 3.4):
\begin{equation}
\label{eqn:end}
\mathrm{dim}\, \mathrm{Hom}(E_{C,A},N_{C|S}\otimes A^*)=h^0(S,F_{C,A}\otimes E_{C,A})= \mathrm{dim}\, \mathrm{Hom}(N\otimes \mathcal I_\xi,M)+1.
\end{equation}

We suppose that $E_{C,A}$ is indecomposable. Fix $N$ and $\ell$ and denote by $\mathcal P_{N,\ell}$ the parameter space of bundles $E$ with Chern classes $c_1(E)=L$ and $c_2(E)=d$ given by extensions of type (\ref{eqn:max dest}), and by $\mathcal G_{N,\ell}$ the Grassmann bundle over $\mathcal P_{N,\ell}$ whose fibre over $E$ is $G(2,H^0(E))$. If we assume that $\mathcal P_{N,\ell}$ contains Lazarsfeld-Mukai bundles corresponding to $\mathcal W$, then we have a rational map
\[
\pi_{N,\ell}:\mathcal{G}_{N,\ell}\dashrightarrow \mathcal W
\]
whose fibre over $E_{C,A}$ is the projectivisation of $\mathrm{Hom}(E_{C,A},N_{C|S}\otimes A^*)$.
Since $\mathrm{Pic}(S)$ is discrete and the Lazarsfeld-Mukai condition is open, it follows that for a given $N$ and $\ell$, the map $\pi_{N,\ell}$ is dominant.
Hence, using (\ref{eqn:end}), it suffices to prove to prove that
\begin{equation}
\label{eqn:G}
\mathrm{dim}\,\mathcal G_{N,\ell}-\mathrm{dim}\, \mathrm{Hom}(N\otimes \mathcal I_\xi,M)\le \mathrm{dim}|L|+(d-k).
\end{equation}

Similarly to \cite[Lemma 3.10]{Aprodu-Farkas}  and \cite[Lemma 4.1]{LelliChiesa}, we obtain the inequality
\begin{equation}
\label{eqn:Cliff}
M\cdot N+\omega_S\cdot N+2\ge k.
\end{equation}

We have
\[
\mathrm{dim}\,\mathcal G_{N,\ell}\le \mathrm{dim}\, G(2,H^0(E_{C,A}))+\mathrm{dim}\, S^{[\ell]}+(\mathrm{dim}\, \mathrm{Ext}^1(N\otimes \mathcal I_\xi,M)-1),
\]
and hence, using Serre duality and observing that $M\cdot H\ge N\cdot H\ge (N+\omega_S)\cdot H\ge 0$ implies $\mathrm{Ext}^2(N\otimes \mathcal I_\xi,M)=\mathrm{Hom}(M,N\otimes\mathcal{I}_\xi\otimes\omega_S)=0$, we obtain the estimate:
\[
\mathrm{dim}\,\mathcal G_{N,\ell}-\mathrm{dim}\, \mathrm{Hom}(N\otimes \mathcal I_\xi,M)\le \mathrm{dim}\, G(2,H^0(E_{C,A}))+2\ell-\chi(S,M^*\otimes N\otimes \omega_S\otimes \mathcal I_\xi)-1.
\]

By the assumption $h^1(S,E_{C,A})=0$, we have $h^0(C,\omega_S\otimes A)=0$, Remark \ref{rmk:pg=0}. It follows that $h^0(C,N_{C|S}\otimes A^*)=g-d-1-\omega_S\cdot C$.
Since $h^0(S,E_{C,A})=2+h^0(C,N_{C|S}\otimes A^*)$ it implies 
\[
\mathrm{dim}\, G(2,H^0(E_{C,A}))=2(g-d-1-\omega_S\cdot C)
\]

The conclude the proof of inequality (\ref{eqn:G}) we compute by the Riemann-Roch theorem (compare to \cite{LelliChiesa})
\[
\chi(S,M^*\otimes N\otimes \omega_S\otimes \mathcal I_\xi)=g-2N\cdot M-\omega_S\cdot M-\ell.
\]
and use (\ref{eqn:Cliff}) and the inequality:
\[
h^0(S,L)\ge g-\omega_S\cdot C=\chi(S,L);
\]
note that $h^2(S,L)=0$, since $L$ is effective on a surface with $p_g=0$.

\medskip

{\em Subcase I.c.} Assume that, for all $(C,A)\in \mathcal W$, $h^1(S,E_{C,A})>0$. From Remark \ref{rmk:pg=0}, we are in the situation $h^0(C,\omega_S\otimes A)>0$, in particular, $A\in \{\omega_S^*|_C\}+W^0_{d+\omega_S\cdot C}(C)$. Since $\omega_S\cdot C\le -k$, it follows that $A$ moves in a family of dimension $\le d-k$.

\medskip

{\em Case II.} Assume that, for any $(C,A)\in \mathcal W$, $A$ has base-points. Then we apply the inductive argument and reduce to the previous case. Note that this case cannot occur for $d=k$.
\endproof

\section{Syzygies of curves}
\label{sec:syz}

In recent curve theory a lot of effort has been put into understanding the relations between syzygies of canonical curves (algebraic objects) and the existence of special linear series (geometric objects). The interest in clarifying these deep relationships between algebraic and geometric properties is high, as failure of vanishing of syzygies produces interesting determinantal cycles on various moduli spaces, and the canonical case is the most natural and basic situation. The precise relationship is predicted by {\em Green's conjecture}: the ideal of a non-hyperelliptic curve is generated by quadrics, and the Clifford dimension controls the number of steps up to which the syzygies are linear. In the language of Koszul cohomology using duality \cite{Green84}, it amounts to the following relation
\[
K_{p,1}(C,K_C)=0\mbox{ for all }p\ge g-c-1,
\]
for any curve $C$ of genus $g$ and Clifford index $c$; for the precise definitions of the objects involved in the statement, we refer to \cite{Green84}. Green's conjecture is known to be true for general curves, \cite{Voisin02,Voisin05}, and moreover the dimension computations of Brill-Noether loci, in particular conditions similar to (\ref{eqn:lgc W}), are related to syzygies of canonical curves. This relationship is explained in \cite{Aprodu05}, \cite{Aprodu13} and has been used in \cite{Aprodu-Farkas} for curves on $K3$ surfaces.
In our case, Green's conjecture is satisfied for {\em general} curves, in the linear system $|L|$, which verify the hypotheses of Theorem \ref{thm:main}. However, under stronger hypotheses, we can prove Green's conjecture for {\em every} curve in the corresponding linear system:

\begin{thm}
\label{thm:Green}
Under the assumptions of Theorem \ref{thm:main}, suppose moreover that 
\[
g-k\ge h^0(S,\omega^{\otimes 2}_S(C)).
\]
Then any smooth curve in $|L|$ is of gonality $k$, Clifford dimension one, and satisfies Green's conjecture.
\end{thm}

\proof
The long exact sequence associated to 
\[
0\to \omega_S\to \omega_S\otimes L\to K_C\to 0
\]
shows that the restriction morphism $H^0(S,\omega_S\otimes L)\to H^0(C,K_C)$ is an isomorphism and provides us with a long exact sequence on Koszul cohomology \cite[(1.d.4)]{Green84}
\[
0= K_{p,1}(S,-C,\omega_S\otimes L)\to K_{p,1}(S,\omega_S\otimes L)\to K_{p,1}(C,K_C)\to K_{p-1,2}(S,-C,\omega_S\otimes L)\to\cdots
\]
Green's vanishing Theorem \cite{Green84} (3.a.1) implies 
\[
K_{p-1,2}(S,-C,\omega_S\otimes L)=0
\]
for $p\ge h^0(S,\omega_S^{\otimes 2}(C))+1$, in particular, for any $C$ and any $p\ge g-k+1$ we obtain an isomorphism
\[
K_{p,1}(S,\omega_S\otimes L)\stackrel{\sim}{\to}K_{p,1}(C,K_C).
\]
Since Green's conjecture is valid for general curves $C\in|L|$ which have gonality $k$ and Clifford dimension one, we infer that $K_{p,1}(S,\omega_S\otimes L)=0$ for any $p\ge g-k+1$.  In particular, $K_{p,1}(C,K_C)=0$ for any $p\ge g-k+1$ and {\em any} smooth curve $C$. From the Green-Lazarsfeld non-vanishing Theorem \cite[Appendix]{Green84}, it follows that {\em any} smooth curve $C$ must have gonality $k$ and Clifford dimension one.
\endproof

\begin{ex}[Smooth curves on Hirzebruch surfaces]
\label{ex:Hirzebruch}
The hypotheses of Theorem \ref{thm:Green} are realised for curves on Hirzebruch surfaces $S=\Sigma_e$ with $e\ge 2$, hence we obtain an alternate proof of the results of~\cite{Aprodu02}.

Indeed, denote by $C_0$ the minimal section and by $F$ the class of a fibre, and let $C\equiv kC_0+mF$ with $m\ge ke$ by a smooth curve on $S$. The gonality of $C$ is $k$ and the genus of $C$ is computed by the formula 
\[
g=(k-1)\left(m-1-\frac{ke}{2}\right).
\]

The first condition $-\omega_S\cdot C\ge k$, from Theorem \ref{thm:main}, is easily verified, as $(2C_0+(e+2)F)\cdot(kC_0+mF)=-ke+2m+2k$. The second condition $g-k\ge h^0(S,\omega_S^{\otimes 2}(C))$, from Theorem \ref{thm:Green}, is verified by direct computation, using the vanishing of $h^1$ of $\omega_S^{\otimes 2}(C)\equiv (k-4)C_0+(m-2e-4)F$ and applying the Riemann-Roch Theorem.
\end{ex}

\end{document}